\documentclass{ifacconf}

\usepackage{graphicx}   
\usepackage{natbib}        
\usepackage{amsmath,amssymb}
\usepackage{algorithm,algorithmic}

\begin{document}
\begin{frontmatter}

\title{Preconditioning for continuation model predictive control
} 


\author[First]{Andrew Knyazev} 
\author[Second]{Alexander Malyshev}

\address[First]{Mitsubishi Electric Research Labs (MERL)
201 Broadway, 8th floor, Cambridge, MA 02139, USA. (e-mail: knyazev@merl.com), \hfill
({http://www.merl.com/people/knyazev}).}
\address[Second]{Mitsubishi Electric Research Labs (MERL)
201 Broadway, 8th floor, Cambridge, MA 02139, USA. (e-mail: malyshev@merl.com)}

\begin{abstract}               
Model predictive control (MPC) anticipates future events to take appropriate control actions.
Nonlinear MPC (NMPC) deals with nonlinear models and/or constraints. A~Continuation/GMRES Method for NMPC, suggested by T. Ohtsuka in 2004, uses the GMRES iterative algorithm to solve a forward difference approximation $Ax=b$ of the original NMPC equations on every time step. We have previously proposed accelerating the GMRES and MINRES convergence by preconditioning the coefficient matrix $A$. We now suggest simplifying the construction of the preconditioner, by approximately solving a forward recursion for the state and a backward recursion for the costate, or simply reusing previously computed solutions.
\end{abstract}

\begin{keyword}
model predictive control, Continuation/GMRES method, preconditioning.
\end{keyword}

\end{frontmatter}

\section{Introduction}

Model Predictive Control (MPC) is an optimal control technology, which is capable to cope with
constrained systems and widely used in industry and academia; see, e.g., \cite{QiBa:03}, \cite{CaBo:04}, and \cite{GrPa:11}. Nonlinear MPC (NMPC) deals with nonlinear models and/or constraints. Main numerical methods applied in NMPC are surveyed by \cite{DiFeHa:09}.

The Continuation/GMRES by \cite{Oht:04} is one of the real-time numerical methods for NMPC. Ohtsuka's method combines several techniques including replacement of inequality constraints by equality constraints, numerical elimination of the state, by the forward recursion, and the costate, by the backward recursion, and the Krylov subspace iterations for solving nonlinear equations via parameter continuation. \cite{TaOh:04} have introduced a preconditioned C/GMRES method, however, their preconditioner
is inefficient.

Our previous work in \cite{KnFuMa:15} extends Ohtsuka's approach in various ways.
The Continuation NMPC (CNPMC) method is formulated for a more general optimal control model
with additional parameters and terminal constraints, which allows us solving minimal time problems.
We also use preconditioners for CNMPC, based on an explicit construction of the Jacobian
matrices at some time steps, improving convergence of the Krylov iterations.
We propose substituting the MINRES iterative solver for GMRES in CNMPC,
reducing the memory requirements and the arithmetic costs per iteration. 

The present note shows how to reduce the cost of the preconditioning setup, by approximating the Jacobian matrix in the Newton iterations. The idea of such an approximation relies on the observation that most entries of the Jacobian weakly depend on small perturbations of the state and costate.
Most columns of the Jacobian can be built from a single instance of the state and costate variables computed, e.g., during generation of the right-hand side of the system solved by the Newton method. Only a small number of columns of the Jacobian, specifically, responsible for treating the terminal constraints and the parameter, is sensitive to changes of the state and costate.
We recalculate the state and costate corresponding just to these sensitive columns.
Moreover, for the purpose of the preconditioner setup, we can, in addition, compute
the state and costate on a coarser grid on the horizon with subsequent linear interpolation of them at the intermediate points. We can also use other general techniques for fast preconditioner setup, e.g.,\ computation of the state and costate variables, as well as the preconditioner and its factorization, in a reduced computer precision. Our numerical results demonstrate that the preconditioned GMRES and MINRES, where the preconditioner is constructed using the approximate state and costate variables, converge faster, compared to their analogs without preconditioning. The paper discusses basic principles of preconditioning,
and detailed algorithms of computation of the preconditioning schemes are to be reported
in our extended paper.

The rest of the note is as follows. In Section~2, we derive the nonlinear equations,
which are solved by the continuation Newton-Krylov method. Section 3 describes how GMRES or MINRES iterations are applied to numerical solution of these nonlinear equations.
Section 4 presents our main contribution by giving details of the preconditioner construction,
which is based on reusing the previously computed and approximated state and costate variables. Section 5 defines a representative test example; and Section 6 gives numerical results illustrating the quality of the method with the suggested preconditioner.

\section{Derivation of the optimality conditions}

The MPC approach is based on the prediction by means of
a finite horizon optimal control problem along a fictitious time $\tau\in[t,t+T]$.
Our model finite horizon problem consists in choosing the control $u(\tau)$
and parameter vector $p$, which minimize the performance index $J$ as follows:
\[
\min_{u,p} J,
\]
where
\[
J = \phi(x(t+T),p)+\int_t^{t+T}L(\tau,x(\tau),u(\tau),p)d\tau
\]
subject to the equation for the state dynamics
\begin{equation}\label{e1}
\frac{dx}{d\tau}=f(\tau,x(\tau),u(\tau),p),
\end{equation}
and the equality constraints for the state $x$ and control $u$
\begin{equation}\label{e2}
C(\tau,x(\tau),u(\tau),p) = 0,
\end{equation}
\begin{equation}\label{e3}
\psi(x(t+T),p) = 0.
\end{equation}
The initial value $x(\tau)|_{\tau=t}$ for (\ref{e1}) is the state vector $x(t)$
of the dynamic system. The control vector $u=u(\tau)$, solving the problem over the prediction
horizon, is used afterwards as an input to control the system at time $t$.
The components of the vector $p(t)$ are parameters of the system and
do not depend on~$\tau$. In our minimum-time example in Section 5,
the scalar parameter $p(t)$ denotes the time to destination, and the horizon
length is $T=p(t)$.

The prediction problem stated above is discretized
on a uniform, for simplicity of presentation, time grid over the horizon $[t,t+T]$
partitioned into $N$ time steps of size $\Delta\tau$, and the time-continuous vector
functions $x(\tau)$ and $u(\tau)$ are replaced by their sampled values $x_i$ and $u_i$
at the grid points $\tau_i$, $i=0,1,\ldots,N$. The integral of the performance cost $J$
over the horizon is approximated by the rectangular quadrature rule.
Equation (\ref{e1}) is integrated by the the explicit Euler scheme,
which is the simplest possible method. We note that more sofisticated
one-step adaptive schemes can be used as well.
The discretized optimal control problem is formulated as follows:
\[
\min_{u_i,p}\left[
\phi(x_N,p) + \sum_{i=0}^{N-1}L(\tau_i,x_i,u_i,p)\Delta\tau\right],
\]
subject to
\begin{equation}\label{e4}
\quad x_{i+1} = x_i + f(\tau_i,x_i,u_i,p)\Delta\tau,\quad i = 0,1,\ldots,N-1,
\end{equation}
\begin{equation}\label{e5}
C(\tau_i,x_i,u_i,p) = 0,\quad  i = 0,1,\ldots,N-1,
\end{equation}
\begin{equation}\label{e6}
\psi(x_N,p) = 0.
\end{equation}

The necessary optimality conditions for the discretized finite horizon problem
are the stationarity conditions for the discrete Lagrangian function
\begin{eqnarray*}
&&\mathcal{L}(X,U)=\phi(x_N,p)+\sum_{i=0}^{N-1}
L(\tau_i,x_i,u_i,p)\Delta\tau\\
&&+\,\lambda_0^T[x(t)-x_0]+\sum_{i=0}^{N-1}\lambda_{i+1}^T[x_i-x_{i+1}+f(\tau_i,x_i,u_i,p)\Delta\tau]\\
&&+\sum_{i=0}^{N-1}\mu_i^TC(\tau_i,x_i,u_i,p)\Delta\tau+\nu^T\psi(x_N,p),
\end{eqnarray*}
where $X = [x_i\; \lambda_i]^T$, $i=0,1,\ldots,N$, and
$U = [u_i\; \mu_i\; \nu\; p]^T$, $i=0,1,\ldots,N-1$.
Here, $\lambda$ is the costate vector, $\mu$ is the Lagrange multiplier vector 
associated with the constraint~(\ref{e5}). The terminal constraint (\ref{e6})
is relaxed by the aid of the Lagrange multiplier $\nu$.

The necessary optimality conditions are the system of nonlinear equations
$\mathcal{L}_{\lambda_i}=0$, $\mathcal{L}_{x_i}=0$, $i=0,1,\ldots,N$,
$\mathcal{L}_{u_j}=0$, $\mathcal{L}_{\mu_j}=0$, $i=0,1,\ldots,N-1$,
$\mathcal{L}_{\nu_k}=0$, $\mathcal{L}_{p_l}=0$.

For further convenience, we introduce the Hamiltonian function
$H(t,x,\lambda,u,\mu,p) = L(t,x,u,p)+\lambda^T f(t,x,u,p)+\mu^T C(t,x,u,p)$.

The optimality conditions are reformulated in terms of a mapping $F[U,x,t]$,
where the vector $U$ combines the control input $u$, the Lagrange multiplier
$\mu$, the Lagrange multiplier $\nu$, and the parameter $p$, all in one vector:
\[
U(t)=[u_0^T,\ldots,u_{N-1}^T,\mu_0^T,\ldots,\mu_{N-1}^T,\nu^T,p^T]^T. 
\]
The vector argument $x$ in $F[U,x,t]$ denotes the state vector at time $t$,
which serves as the initial vector $x_0$ in the following procedure.
\begin{enumerate}
\item Starting from the current measured or estimated state $x_0$, compute all
$x_i$, $i=0,1\ldots,N$, by the forward recursion
\[
x_{i+1} = x_i + f(\tau_i,x_i,u_i,p)\Delta\tau.
\]
Then starting from
\[
\lambda_N=\frac{\partial\phi^T}{\partial x}(x_N,p)+
 \frac{\partial\psi^T}{\partial x}(x_N,p)\nu
\]
compute all costates $\lambda_i$, $i=N,\ldots,1,0$, by the backward recursion
\[
\lambda_i=\lambda_{i+1}+\frac{\partial H^T}{\partial x}
(\tau_i,x_i,\lambda_{i+1},u_i,\mu_i,p)\Delta\tau.
\]
\item Using just obtained $x_i$ and $\lambda_i$, calculate the vector
\begin{eqnarray*}
&&\hspace*{-2.6em}F[U,x,t]=\left[\begin{array}{c}\begin{array}{c}
\frac{\partial H^T}{\partial u}(\tau_0,x_0,\lambda_{1},u_0,\mu_0,p)\Delta\tau\\
\vdots\\\frac{\partial H^T}{\partial u}(\tau_i,x_i,\lambda_{i+1},u_i,\mu_i,p)\Delta\tau\\
\vdots\\\frac{\partial H^T}{\partial u}(\tau_{N-1},x_{N-1},\lambda_{N},u_{N-1},
\mu_{N-1},p)\Delta\tau\end{array}\\\;\\
\begin{array}{c}C(\tau_0,x_0,u_0,p)\Delta\tau\\
\vdots\\C(\tau_i,x_i,u_i,p)\Delta\tau\\\vdots\\
C(\tau_{N-1},x_{N-1},u_{N-1},p)\Delta\tau\end{array}\\\;\\
\psi(x_N,p)\\[2ex]
\begin{array}{c}\frac{\partial\phi^T}{\partial p}(x_N,p)+
\frac{\partial\psi^T}{\partial p}(x_N,p)\nu\\
+\sum_{i=0}^{N-1}\frac{\partial H^T}{\partial p}(\tau_i,x_i,
\lambda_{i+1},u_i,\mu_i,p)\Delta\tau\end{array}
\end{array}\right]\!\!.
\end{eqnarray*}
\end{enumerate}
The equation with respect to the unknown vector $U(t)$
\begin{equation}\label{e7}
 F[U(t),x(t),t]=0
\end{equation}
gives the required necessary optimality conditions.

\section{Numerical algorithm}
The controlled system is sampled on a uniform time grid $t_j=j\Delta t$, $j=0,1,\ldots$.
Solution of equation (\ref{e7}) must be found at each time step $t_j$ in real time,
which is a challenging part of implementation of NMPC.

Let us denote $x_j=x(t_j)$, $U_j=U(t_j)$, and rewrite the equation $F[U_j,x_j,t_j]=0$
equivalently in the form
\[
F[U_j,x_j,t]-F[U_{j-1},x_j,t_j]=b_j,
\]
where
\begin{equation}\label{e8}
b_j=-F[U_{j-1},x_j,t_j].
\end{equation}

Using a small $h$, which may be different from $\Delta t$ and $\Delta\tau$,
we introduce the operator
\begin{eqnarray}\label{e9}
a_j(V)=(F[U_{j-1}+hV,x_j,t_j]-F[U_{j-1},x_j,t_j])/h.
\end{eqnarray}
We note that the equation $F[U_j,x_j,t_j]=0$ is equivalent to the equation
$a_j(\Delta U_j/h)=b_j/h$, where $\Delta U_j=U_j-U_{j-1}$.

Let us denote the $k$-th column of the $m\times m$ identity matrix by $e_k$,
where $m$ is the dimension of the vector $U$, and define an $m\times m$ matrix $A_j$
with the columns $A_je_k=a_j(e_k)$, $k=1,\ldots,m$.
The matrix $A_j$ is an $O(h)$ approximation of the Jacobian matrix
$F_U[U_{j-1},x_j,t_j]$, which is symmetric.

Suppose that an approximate solution $U_0$ to the initial equation $F[U_0,x_0,t_0]=0$ is
available. The first block entry of $U_0$ is taken as the control $u_0$
at the state $x_0$. The next state $x_1=x(t_1)$ is either sensor estimated or computed
by the formula $x_1=x_0+f(t_0,x_0,u_0)\Delta t$; cf. (\ref{e1}).

At the time $t_j$, $j>1$, we have the state $x_j$ and the vector $U_{j-1}$
from the previous time $t_{j-1}$. Our goal is to solve the following equation
with respect to $V$:
\begin{equation}\label{e11}
 a_j(V)=b_j/h.
\end{equation}
Then we can set $U_j=U_{j-1}+hV$ and choose the first block component of $U_j$ as the control $u_j$.
The next system state $x_{j+1}=x(t_{j+1})$
is either sensor estimated or computed by the formula $x_{j+1}=x_j+f(t_j,x_j,u_j)\Delta t$.

A direct way to solve (\ref{e11}) is generating the matrix $A_j$
and then solving the system of linear equations $A_j\Delta U_j=b_j$;
e.g., by the Gaussian elimination.

A less expensive alternative is solving (\ref{e11}) by the GMRES method,
where the operator $a_j(V)$ is used without explicit construction of
the matrix $A_j$ (cf., \cite{Kel:95,Oht:04}). Some results on convergence
of GMRES in the nonlinear case can be found in \cite{BrWaWo:08}.

We recall that, for a given system of linear equations $Ax=b$ and initial approximation $x_0$,
GMRES constructs orthonormal bases of the Krylov subspaces
$\mathcal{K}_n=\text{span}\{r_0,Ar_0,\ldots,A^{n-1}r_0\}$,
$n=1,2,\ldots$, given by the columns of matrices $Q_n$,
such that $AQ_n=Q_{n+1}H_{n}$ with the upper Hessenberg matrices~$H_n$
and then searches for approximations to the solution $x$ in the form $x_n=Q_ny_n$,
where $y_n=\text{argmin}\|AQ_ny_n-b\|_2$.

A more efficient variant of GMRES, called MINRES, may be applied when the matrix $A$
is symmetric, and the preconditioner is symmetric positive definite. 
Using the MINRES iteration in Ohtsuka's approach is mentioned in \cite{KnFuMa:15}.

\section{Preconditioning}

The convergence of GMRES can be accelerated by preconditioning. A matrix $M$ that is close
to the matrix $A$ and such that computing $M^{-1}r$ for an arbitrary vector $r$
is relatively easy, is referred to as a preconditioner. The preconditioning
for the system of linear equations $Ax=b$ with the preconditioner $M$ formally
replaces the original system $Ax=b$ with the equivalent preconditioned linear system
$M^{-1}Ax=M^{-1}b$. If the condition number $\|M^{-1}A\|\|A^{-1}M\|$
of the matrix $M^{-1}A$ is small, convergence of iterative Krylov-based solvers
for the preconditioned system can be fast. However, in general, the convergence speed of, e.g.,\ the preconditioned GMRES is not necessarily determined by the condition number alone. 

A typical implementation of the preconditioned GMRES is given below.
The unpreconditioned GMRES is the same algorithm but with $M=I$, where $I$ is the identity matrix.
We denote by $H_{i_1:i_2,j_1:j_2}$ the submatrix of $H$ with the entries $H_{ij}$
such that $i_1\leq i\leq i_2$ and $j_1\leq j\leq j_2$.

\begin{tabular}{l}
\hline\\[-2ex]
\textbf{Algorithm}{ Preconditioned GMRES($k_{\max}$)}\\\hline\\[-2ex]
\textbf{Input:} $a(v)$, $b$, $x_0$, $k_{\max}$, $M$\\
\textbf{Output:} Solution $x$ of $a(x)=b$\\
\quad $r=b-a(x_0)$, $z=M^{-1}r$, $\beta=\|z\|_2$, $v_1=z/\beta$\\
\quad\textbf{for }{$k=1,\ldots,k_{\max}$}\textbf{ do}\\
\qquad $r=a(v_k)$, $z=M^{-1}r$\\
\qquad $H_{1:k,k}=[v_1,\ldots,v_k]^Tz$\\
\qquad $z=z-[v_1,\ldots,v_k]H_{1:k,k}$\\
\qquad $H_{k+1,k}=\|z\|_2$\\
\qquad $v_{k+1}=z/\|z\|_2$\\
\quad\textbf{end for}\\
\quad $y=\mbox{arg min}_y\|H_{1:k_{\max}+1,1:k_{\max}}y-[\beta,0,\dots,0]^T\|_2$\\
\quad $x=x_0+[v_1,\ldots,v_{k_{\max}}]y$\\
\end{tabular}

In \cite{KnFuMa:15}, the matrix $A_j$ is exactly computed at some time instances $t_j$
and used as a preconditioner in a number of subsequent time
instances $t_j$, $t_{j+1}$, \ldots, $t_{j+j_p}$.
In the present note, we propose to use a close approximation to $A_j$, which needs
much less arithmetic operations for its setup. Construction of such approximations $M_j$
is the main result of this note.

We recall that computation of the $k$-th column of $A_j$ requires computation
of all states $x(\tau_i)$ and costates $\lambda(\tau_i)$ for the parameters
stored in the vector $U_{j-1}+he_k$. Is it possible to replace them by
$x(\tau_i)$ and $\lambda(\tau_i)$ computed for the parameters stored
in the vector $U_{j-1}$?
The answer is yes, for the indices $k=1, \ldots, m-l$, where $l$ is the sum of dimensions
of $\psi$ and $p$. These $k$ indices correspond to the terms containing the factor $\Delta\tau$ in the Lagrangian $\mathcal{L}$.

The first $m-l$ columns (and rows, since the preconditioner $M_j$ is symmetric) 
are calculated by the same formulas as those in $A_j$, but with the values $x(\tau_i)$ and
$\lambda(\tau_i)$ computed only once for the parameters stored in the vector $U_{j-1}$,
i.e., when computing the vector $b_j$.
Thus, the setup of $M_j$ computes the states $x(\tau_i)$ and costates
$\lambda(\tau_i)$ only $l$ times instead of $m$ times as for the matrix $A_j$.
It is this reduction of computing time that makes the preconditioner $M_j$ more 
efficient, especially in cases where dimension of the state space is very large.

The preconditioner $M_j$ is obtained from $A_j$ by neglecting the derivatives
$\partial x_{i_1}/\partial u_{i_2}$, $\partial\lambda_{i_1}/\partial u_{i_2}$,
$\partial x_{i_1}/\partial\mu_{i_2}$ and $\partial\lambda_{i_1}/\partial\mu_{i_2}$.
Therefore, the difference $A_j-M_j$ is of order $O(\Delta\tau)$ since 
$\partial x_{i_1}/\partial u_{i_2}=O(\Delta\tau)$,
$\partial\lambda_{i_1}/\partial u_{i_2}=O(\Delta\tau)$,
$\partial x_{i_1}/\partial\mu_{i_2}=0$ and
$\partial\lambda_{i_1}/\partial\mu_{i_2}=O(\Delta\tau)$.

The preconditioner application $M^{-1}r$ requires the LU factorization $M=LU$, which
is computed by the Gaussian elimination. Then the vector $M^{-1}r=U^{-1}(L^{-1}r)$ is
obtained by performing back-substitutions for the triangular factors $L$ and $U$.
Further acceleration of the preconditioner setup is possible by faster computation of
the LU factorization. For example, when computation with lower number of bits is cheaper
than computation with the standard precision, the preconditioner $M_j$
and its LU factorization may be computed in lower precision.

Another way of reduction of the arithmetical work in the preconditioner setup
is the computation of the states $x$ and costates $\lambda$ with the double step
$2\Delta\tau$ thus halving the arithmetical cost and memory storage.
The intermediate values of $x$ and $\lambda$ are then obtained from the computed
values by simple linear interpolation.

\section{Example}

We consider a test nonlinear problem, which describes the minimum-time motion from a state
$(x_0,y_0)$ to a state $(x_f,y_f)$ with an inequality constrained control:
\begin{itemize}
\item State vector
 $\vec{x}=\left[\begin{array}{c}x\\y\end{array}\right]$ and
input control $\vec{u}=\left[\begin{array}{c}u\\u_d\end{array}\right]$.
\item Parameter variable $\vec{p}=t_f-t$, where $t_f$ denotes the arrival time
at the terminal state $(x_f,y_f)$.
\item Nonlinear dynamics is governed by the system of ordinary differential equations
\[
\dot{\vec{x}}=f(\vec{x},\vec{u},\vec{p})=
\left[\begin{array}{c}(Ax+B)\cos u\\(Ax+B)\sin u\end{array}\right].
\]
\item Constraint: $C(\vec{x},\vec{u},\vec{p})=(u-c_{u})^2+u_d^2-r_{u}^2=0$,
where $c_u=c_0+c_1\sin(\omega t)$ and $u_d$ is a slack variable,
i.e., the control $u$ always stays within
the sinusoidal band $c_{u}-r_{u}\leq u\leq c_{u}+r_{u}$).
\item Terminal constraints: $\psi(\vec{x},\vec{p})=\left[\begin{array}{c}
x-x_f\\y-y_f\end{array}\right]=0$ (the state should pass through the point
$(x_f,y_f)$ at $t=t_f$)
\item Objective function on the horizon interval $[t,t_f]$:
\[
J=\phi(\vec{x},\vec{p})+\int_t^{t+\vec{p}}L(\vec{x},\vec{u},\vec{p})dt,
\]
where
\[
\phi(\vec{x},\vec{p})=\vec{p},\quad L(\vec{x},\vec{u},\vec{p})=-w_{d}u_d
\]
(the state should arrive at $(x_f,y_f)$ in the shortest time;
the function $L$ serves to stabilize the slack variable $u_d$)
\item Constants: $A=B=1$, $x_0=y_0=0$, $x_f=y_f=1$,
$c_0=0.8$, $c_1=0.3$, $\omega=10$, $r_{u}=0.2$, $w_{d}=0.005$.
\end{itemize}
The horizon interval $[t,t_f]$ is parametrized by the affine mapping
$\tau\to t+\tau\vec{p}$ with $\tau\in[0,1]$.

The components of the corresponding discretized problem on the horizon are given below:
\begin{itemize}
\item $\Delta\tau=1/N$, $\tau_i=i\Delta\tau$,
$c_{ui}=c_0+c_1\sin(\omega(t+\tau_ip))$;
\item the participating variables are the state $\left[\begin{array}{c}
x_i\\y_i\end{array}\right]$, the costate $\left[\begin{array}{c}
\lambda_{1,i}\\\lambda_{2,i}\end{array}\right]$, the control $\left[\begin{array}{c}
u_{i}\\u_{di}\end{array}\right]$, the Lagrange multipliers
$\mu_i$ and $\left[\begin{array}{c}\nu_{1}\\\nu_{2}\end{array}\right]$,
the parameter $p$;
\item the state is governed by the model equation
\[
\left\{\begin{array}{l} x_{i+1}=x_i+\Delta\tau\left[p\left(Ax_{i}+B\right)\cos u_{i}\right],\\
y_{i+1}=y_i+\Delta\tau\left[p\left(Ax_{i}+B\right)\sin u_{i}\right],\end{array}\right.
\]
where $i=0,1,\ldots,N-1$;
\item the costate is determined by the backward recursion ($\lambda_{1,N}=\nu_1$,
$\lambda_{2,N}=\nu_2$)
\[
\left\{\begin{array}{l} \lambda_{1,i}=\lambda_{1,i+1}
\\\hspace{2.5em}{}
+\Delta\tau\left[pA(\cos u_i \lambda_{1,i+1}+\sin u_i\lambda_{2,i+1})\right],\\
\lambda_{2,i} = \lambda_{2,i+1},\end{array}\right.
\]
where $i=N-1,N-2,\ldots,0$;
\item the equation $F(U,x_0,y_0,t)=0$, where
\begin{eqnarray*}
\lefteqn{U=[u_0,u_{d,0},\ldots,u_{N-1},u_{d,N-1},}
&\hspace*{10em}\\&&\mu_0,\ldots,\mu_{N-1},\nu_1,\nu_2,p],
\end{eqnarray*}
has the following rows from the top to bottom:
\[
\left\{\begin{array}{l}
\Delta\tau\left[p(Ax_i+B)\left(-\sin u_i\lambda_{1,i+1}+
\cos u_i\lambda_{2,i+1}\right)\right.\\
\hspace*{11em}\left.{}+2\left(u_i-c_{ui}\right)\mu_i\right] = 0 \\
\Delta\tau\left[2\mu_iu_{di}-w_{d}p\right] = 0 \end{array}\right.
\]
\[
\left\{\;\Delta\tau\left[(u_i-c_{ui})^{2}+u_{di}^2-r_{u}^2\right]=0\right.\hspace*{7em}
\]
\[
\left\{\begin{array}{l} x_N-x_r=0\\y_N-y_r=0\end{array}\right.\hspace{15em}
\]
\[
\left\{\begin{array}{l}\Delta\tau[\sum\limits^{N-1}_{i=0}
(Ax_i+B)(\cos u_i\lambda_{1,i+1}+\sin u_i\lambda_{2,i+1})\\
\hspace{3em}{}-2(u_i-c_{ui})\mu_ic_1\cos(\omega(t+\tau_ip))\omega\tau_i\\
\hspace{13em}-w_du_{di}]+1 = 0.\end{array}\right.
\]
\end{itemize}

Let us compare the computation costs of the matrices $A_j$ and $M_j$ for this example.
We do not take into account the computation of the right-hand side $b_j=-F[U_{j-1},x_j,t_j]$
because it is a necessary cost. Computation of the matrix $A_j$ requires $3N+3$ evaluations
of the vector $F[U_{j-1}+hV,x_j,t_j]$, where $N$ is the number of grid points on the prediction horizon.
Setup of $M_j$ requires only 3 evaluations of $F[U_{j-1}+hV,x_j,t_j]$, which is $N+1$ times faster. 

\section{Numerical results}

In our numerical experiments, the weakly nonlinear system~(\ref{e11}) for the test problem
from Section 5 is solved by the GMRES and MINRES iterations.
The number of evaluations of the vector $a(V)$ at each time $t_j$ does not
exceed an a priori chosen constant $k_{\max}=20$. In other words, the maximum number of
GMRES or MINRES iterations is less or equal $k_{\max}$.
The error tolerance in GMRES and MINRES is $tol=10^{-5}$.
The number of grid points on the horizon is $N=50$, the sampling time
of simulation is $\Delta t=1/500$, and $h=10^{-8}$. 

The preconditioners are set up at the time instances $lt_p$, where $t_p=0.2$ is
the period, and $l=0,1,\ldots$. After each setup, the same preconditioner is applied
until next setup. Preconditioners for MINRES must be symmetric positive definite and are built here
as the absolute value of $M_j$, i.e., if $M_j=U\Sigma V^T$ is the singular value decomposition, then
$|M_j|=U\Sigma U^T$; see \cite{VeKn:13}.

\begin{figure}[ht]
\noindent\centering{
\includegraphics[width=3.in]{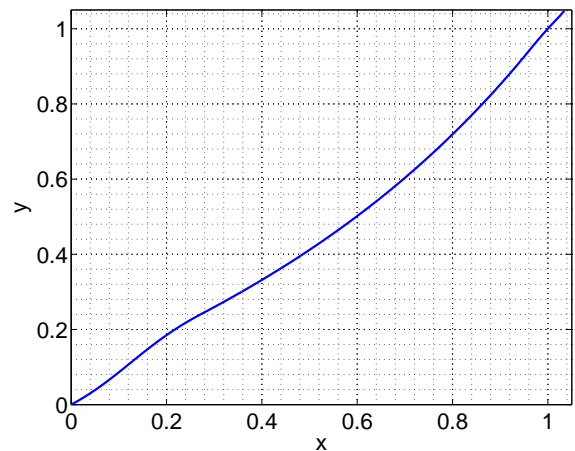}
}
\caption{Trajectory by NPMC using the preconditioned GMRES}
\label{fig1}
\end{figure}

\begin{figure}[ht]
\noindent\centering{
\includegraphics[width=3.in]{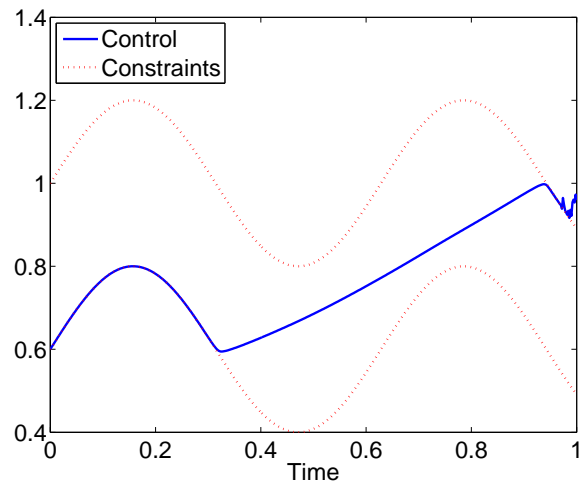}
}
\caption{NMPC control $u$ using the preconditioned GMRES}
\label{fig2}
\end{figure}

\begin{figure}[ht]
\noindent\centering{
\includegraphics[width=3.in]{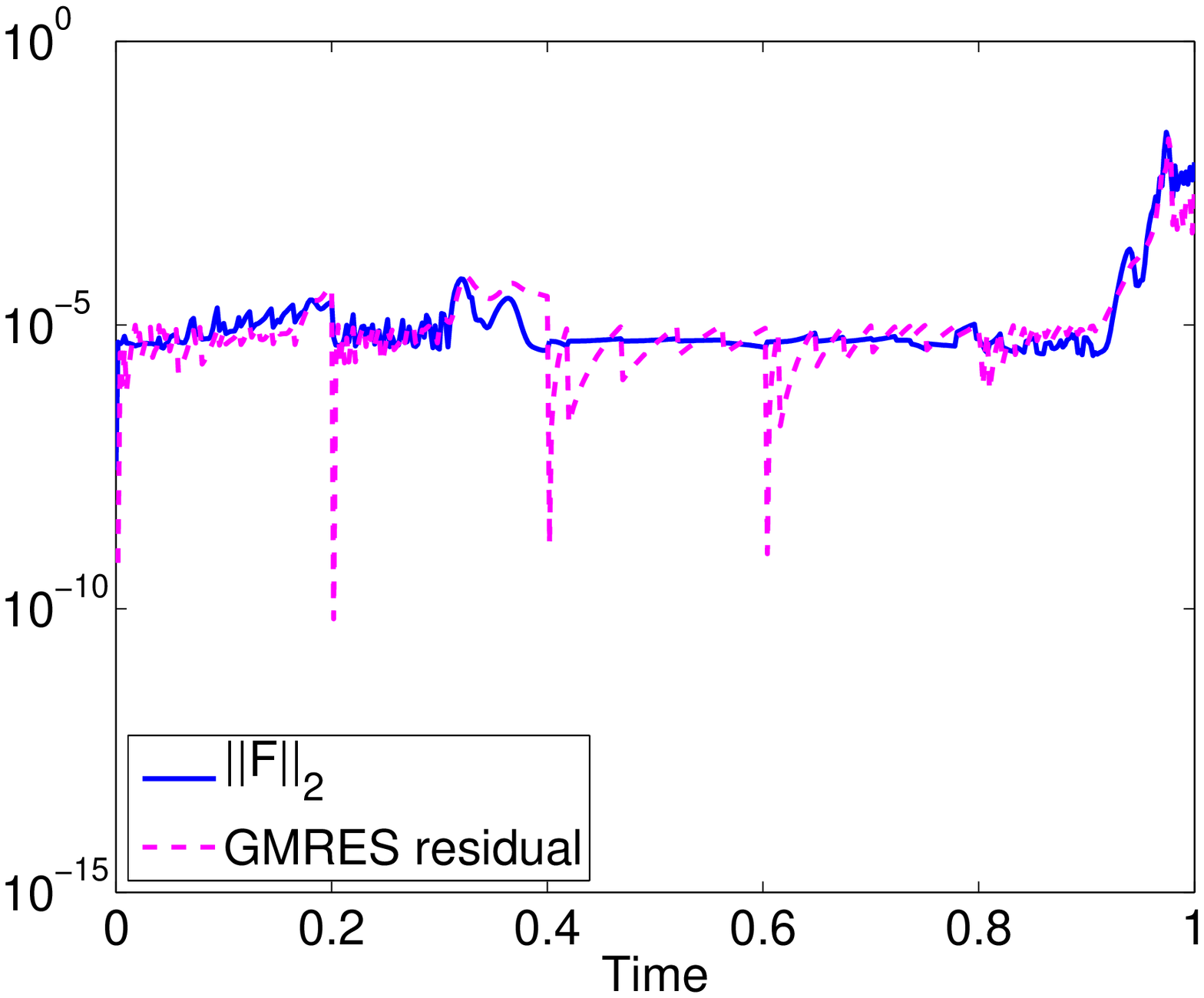}
}
\caption{Preconditioned GMRES, $k_{\max}=20$}
\label{fig3}
\end{figure}

Figure \ref{fig1} shows the computed trajectory for the test example.
Figure \ref{fig2} shows the optimal control by the MPC approach using
the preconditioned GMRES. Figure~\ref{fig3} displays $\|F\|_2$ and the GMRES residuals.

\begin{figure}[ht]
\noindent\centering{
\includegraphics[width=3.in]{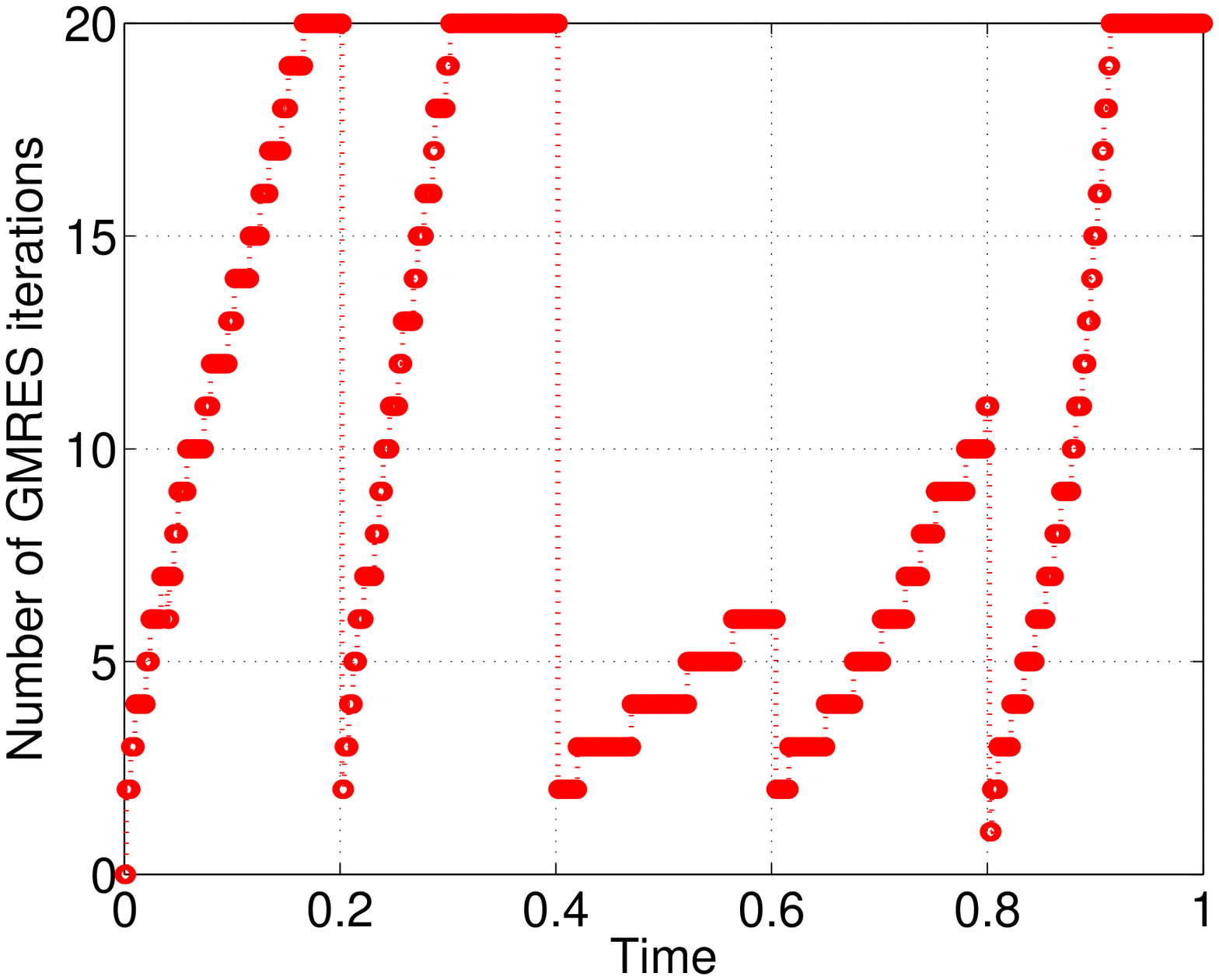}
}
\caption{Preconditioned GMRES, $k_{\max}=20$}
\label{fig4}
\end{figure}

\begin{figure}[ht]
\noindent\centering{
\includegraphics[width=3.in]{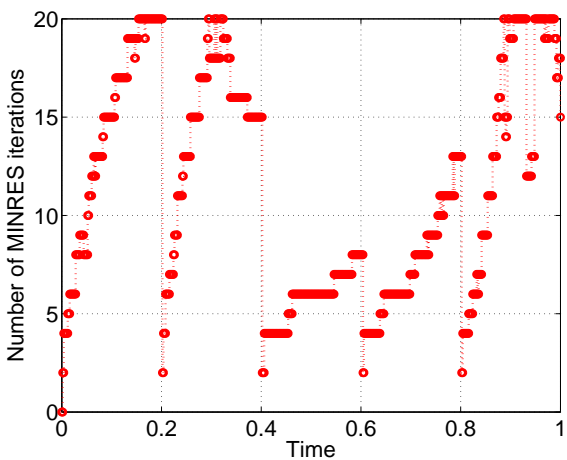}
}
\caption{Preconditioned MINRES, $k_{\max}=20$}
\label{fig5}
\end{figure}

The number of iterations of preconditioned GMRES is displayed
in Figure \ref{fig4}. For comparison, we show the number of iterations
of the preconditioned MINRES in Figure \ref{fig5}. 

\begin{figure}[ht]
\noindent\centering{
\includegraphics[width=3.in]{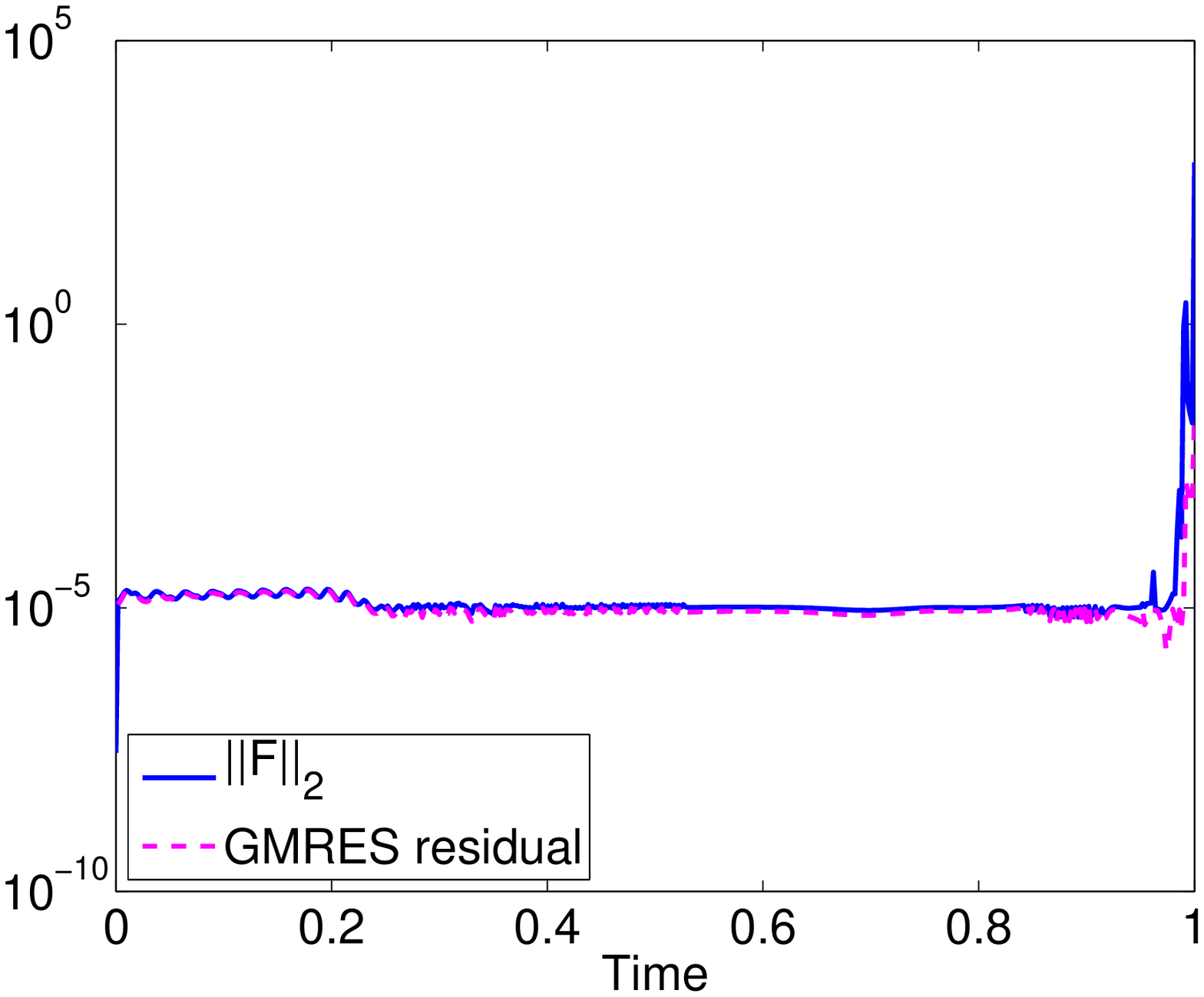}
}
\caption{GMRES without preconditioning, $k_{\max}=20$}
\label{fig6}
\end{figure}

\begin{figure}[ht]
\noindent\centering{
\includegraphics[width=3.in]{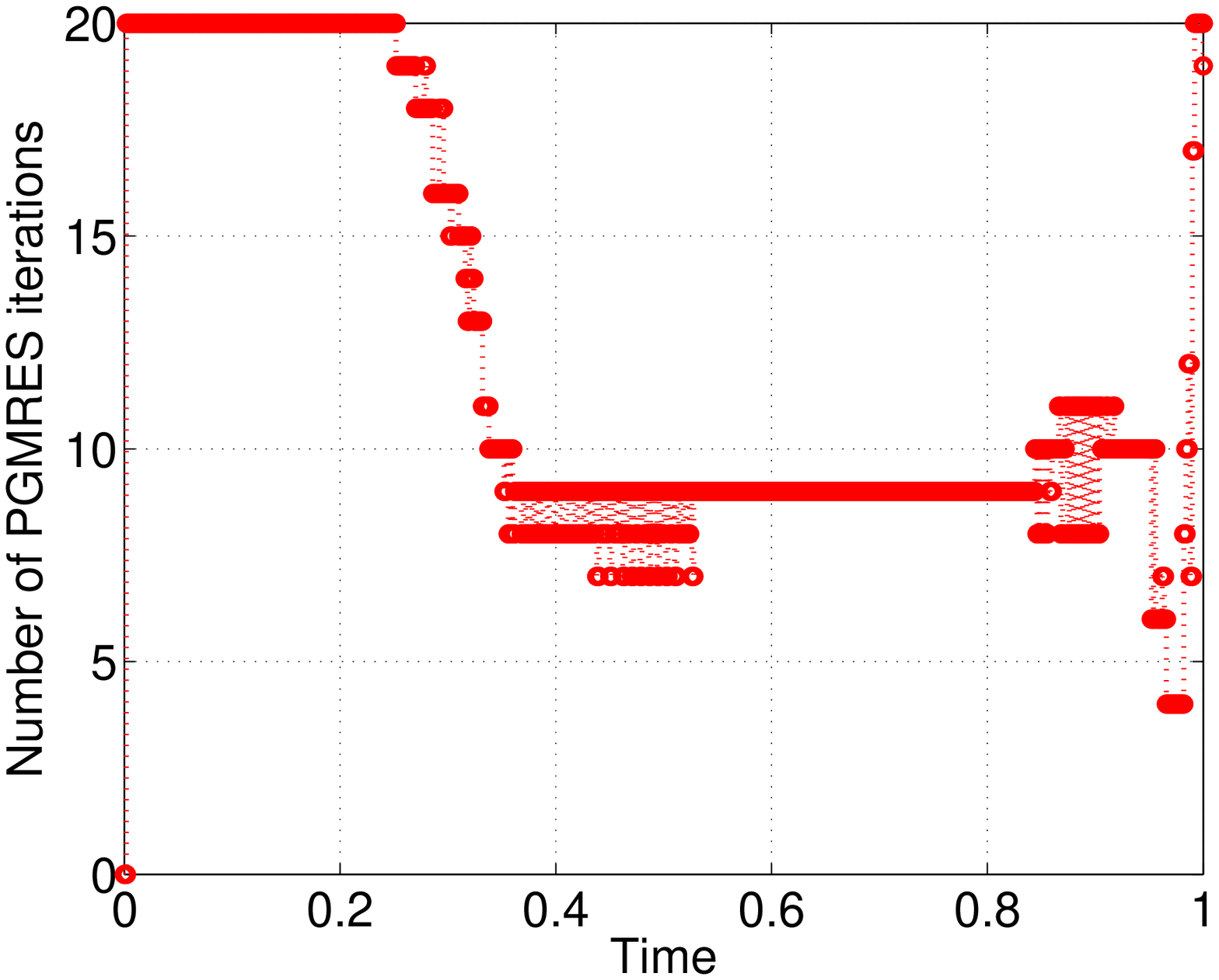}
}
\caption{GMRES without preconditioning, $k_{\max}=20$}
\label{fig7}
\end{figure}

\begin{figure}[ht]
\noindent\centering{
\includegraphics[width=3.in]{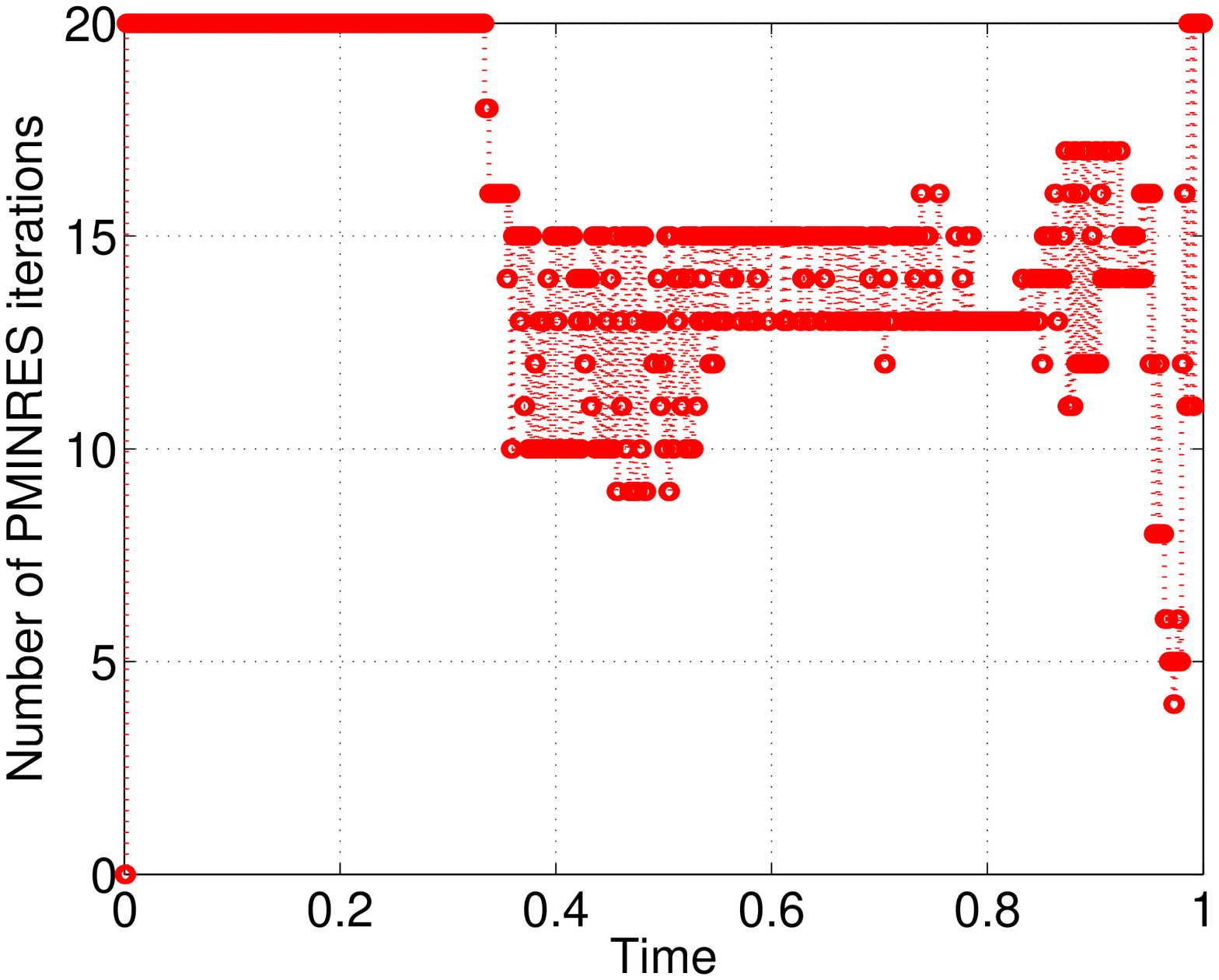}
}
\caption{MINRES without preconditioning, $k_{\max}=20$}
\label{fig8}
\end{figure}

Figure \ref{fig6} displays $\|F\|_2$ and the 2-norm of the residual after iterations
of GMRES without preconditioning. The corresponding number of iterations of GMRES
without preconditioning is shown in Figure \ref{fig7}. The number of iterations of MINRES
without preconditioning is shown in Figure \ref{fig8}.

Effect of preconditioning is seen when comparing Figures~4 and 7 for GMRES and
Figures 5 and 8 for MINRES. Preconditioning of GMRES reduces the number of iterations
by factor 1.2. Preconditioning of MINRES reduces the number of iterations by factor 1.4.

The number of iterations does not necessarily account for the additional complexity
that preconditioning brings to the on-line algorithm. However, the computation time
is machine and implementation dependent, while our tests are done in MATLAB
on a generic computer. Specific implementations on dedicated computer chips for
on-line controllers is a topic of future work. 

\section*{Conclusions}

We have found a new efficient preconditioner $M_j$, which approximates the Jacobian matrix
$A_j$ of the mapping $F$ defining equation (\ref{e7}). Computation of $M_j$ is $O(N)$ times
faster than that of $A_j$, where $N$ is the number of grid points on the prediction horizon.

The preconditioner $M_j$ can be very efficient for the NMPC problems,
where dimension of the state space is large, for example,
in the control of dynamic systems described by partial differential equations.

Other useful techniques for accelerating the preconditioner setup include
computation of the matrices $M_j$ and their LU factorizations
in lower precision, computation of the state and costate on a coarse grid
over the horizon and linear interpolation of the computed values on the fine grid.

\end{document}